
\input amstex
\documentstyle{amsppt}
\pageno=1

\topmatter
 
\title
Equivariant $D$-modules
\endtitle
 
\author Ryoshi Hotta \endauthor
\affil Mathematical Institute, Faculty of Science\\
Tohoku University, Sendai 980, Japan\\
\\
\\
Dedicated to Professor Takeshi Kotake\\
on his sixtieth birthday
\endaffil
 
\abstract
The first part of these notes is devoted to an introduction to algebraic
$D$-modules. Several basic notions as in [2], [11] are introduced.
In the second part, $D$-modules with group action are treated.
Several important examples in this situation are discussed in details.
Particularly, the Harish-Chandra systems for group characters and the
Gelfand generalized hypergeometric systems are our main topics.
\endabstract
 
\thanks
Partly supported by the Grants-in-Aid for Scientific as well as
Co-operative Research, The Ministry of Education, Science and Culture,
Japan.
\endthanks
 
\endtopmatter
 
\def\CC{\Bbb C}
\def\OO{\Cal O}

\document
 
\vskip 10mm
 
\centerline{\bf I. $D$-modules, an introduction}
 
\vskip 10mm
\noindent
\bf 1. Systems of linear partial differential equations \rm
 
\bigskip
 
Let $U$ be a complex domain in the $n$-dimensional complex affine space
${\Bbb C}^n$ and $D(U)$ the ring of partial differential operators
on $U$ with holomorphic coefficients. Consider a system of linear
partial differential equations
 
$$ P_i \, u=0  \quad (1 \leq i \leq m) $$
for $P_i \in D(U)$.
 
Let $F$ be a suitable function space on $U$ stable by the action of
$D(U)$, e.g., $\Cal O(U)$ the space of holomorphic functions,
$C^{\infty}(U)$ that of $C^{\infty}$ functions or $\Cal D'(U)$
that of Schwarz distributions. If $\phi \in F$ is a solution to the above
system of equations $( P_i \, \phi = 0 \  (1 \leq i \leq m))$,
then the map
$$
\tilde \phi : D(U) \ni Q \longmapsto Q \phi \in F
$$
is a left $D(U)$-linear by definition and $\roman{Ker} \, \tilde \phi$
contains the $P_i$'s $(1 \leq i \leq m)$.
Thus the $D(U)$-homomorphism $\tilde \phi$ factorizes to
the $D(U)$-homomorphism
$$
\tilde \phi : D(U)/I \longrightarrow F \quad (Q\mod I \mapsto Q \phi)
$$
where $I = \sum_{i=1}^{m} D(U)P_i$ is the left ideal of the ring $D(U)$
generated by the $P_i$'s.
 
Thus if we denote by $M$ the left $D(U)$-module $D(U)/I$, the space of
solutions to the system in $F$ is identified with the space of
left $D(U)$-module homomorphisms
$$ \roman{Hom}_{D(U)}(M, F) $$
by the correspondence $\phi \leftrightarrow \tilde \phi$.
 
There are several reasons why we consider such algebraic objects,
$D$-modules. First of all, an interpretation of solution spaces
as $\roman{Hom}_D(\ , \ )$ prolongs naturally to use of homological
algebra, which benefits us much enough. Secondly, as will be noted later,
one of the basic invariants, the characteristic variety of a system
can be correctly defined only when we consider the ideal generated
by the $P_i$'s (a fixed set of generators is not enough for the
definition).
 
\vskip 10mm
 
\noindent
{\bf 2. Algebraic differential operators}
\bigskip
 
Since all substantial examples in these notes are algebraic $D$-modules,
we begin with basic notions on algebraic differential operators.
 
Simplest but important examples are linear differential operators with
polynomial coefficients. The ring of differential operators with
polynomial coefficients on the $n$-dimensional complex affine space
${\Bbb C}^n$, denoted by $D({\Bbb C}^n)$, is called the Weyl algebra.
The Weyl algebra $D({\Bbb C}^n)$ is a ${\Bbb C}$-algebra generated by
$$
x_i, \quad \partial_i = \frac{\partial}{\partial x_i} \quad
(1 \leq i \leq n)
$$
with Heisenberg commutator relations
$$
[\partial_i, x_j] = \delta_{ij}, \quad
[x_i, x_j] = [\partial_i, \partial_j] = 0 .
$$
 
Even on general smooth algebraic varieties, the situation does not differ
much from the above. Let $X$ be a smooth affine algebraic variety over
${\Bbb C}$. This means the following. Let $A$ be a commutative algebra
finitely genrated over ${\Bbb C}$ with no nilpotent elements.
The smoothness means that $\dim_{\Bbb C}{\frak m}/{\frak m}^2$ is constant
($= \dim X$) for every maximal ideal $\frak m$ of $A$. The space $X$
is identified with
$\roman{Hom}_{{\Bbb C}-\roman{alg}}(A, {\Bbb C})$, the set
of all ${\Bbb C}$-algebra homomorphisms, which is also identified with
$\roman{Specm} \, A$, the set of all maximal ideals of A by Hilbert's
Nullstellensatz ($ x \leftrightarrow \roman{Ker} \, x = {\frak m}_x
(x \in \roman{Hom}_{{\Bbb C}-\roman{alg}}(A, {\Bbb C})$).
The ${\Bbb C}$-algebra $A$ is then denoted by ${\Bbb C}[X]$ and
called the algebra of regular functions on $X$ ($ f(x) = x(f)$ for
$f \in {\Bbb C}[X], x \in
\roman{Hom}_{{\Bbb C}-\roman{alg}}({\Bbb C}[X], {\Bbb C})$).
The family of subsets $X_f = \{ x \in X | f(x) \neq 0 \}
\  (f \in {\Bbb C}[X])$ forms a basis of open sets in $X$
(the Zariski topology of $X$). Note that ${\Bbb C}[X_f] = {\Bbb C}[X]_f
= {\Bbb C}[X][f^{-1}]$ is the algebra of regular functions of an open
affine subvariety $X_f$ of $X$.
 
The correspondence
$$ X_f \longmapsto \CC[X_f] $$
gives rise to the structure sheaf ${\Cal O}_X$ of $X$ as a local ringed
space (${\Cal O}_X(X_f) \break
= \Gamma (X_f, {\Cal O}_X) = \CC [X_f]$).
The stalk ${\Cal O}_{X,x}$ of ${\Cal O}_X$ at $x \in X$ is the localization
of $\CC[X]$ at the maximal ideal ${\frak m}_x \in \roman{Specm} \, \CC[X] \  
(\displaystyle{\varinjlim_{x \in X_f} \CC[X_f]} = {\Cal O}_{X,{\frak m}_x}))$.
 
In general, a smooth algebraic variety is defined to be a local ringed space
$(X, {\Cal O}_X)$ such that every $x \in X$ has an open neighborhood $U$
such that $(U, {\Cal O}_X|_U)$ is isomorphic to a smooth affine variety
as local ringed spaces as above. (Usually one adopts a further
assumption, i.e.,
separability of the Zariski topology, which means that the diagonal map
$X \overset \Delta\to\rightarrow X \times X \  (\Delta(x) = (x, x))$ is a
closed immersion.)
 
Linear differential operators are defined as follows
in algebraic geometry.
 
\noindent
{\bf Definition.} A $\CC$-linear sheaf endomorphism $P \in
\roman{End}_{\CC}{\Cal O}_X$ is called a {\it linear differential operator of
order not greater than} $m$ if
$$ (\roman{ad} {\Cal O}_X)^{m+1}P = 0. $$
More precisely, for every open $U \subset X$, $P$ is a collection of
${\CC}$-linear maps
$$ P_U \in \roman{End}_{\CC}{\Cal O}_X(U) $$
compatible with all sheaf restriction data ${\Cal O}_X(U) \rightarrow
{\Cal O}_X(V) \  (V \subset U)$ satisfying
$$
[f_0, [f_1, [ \cdots ,[f_m, P_U] \cdots ] = 0 \quad
\text{for every} \  f_0, f_1, \cdots, f_m \in {\Cal O}_X(U).
$$
 
By definition, if  $X$ is affine, a linear differential operator $P$
of order not greater than $m$ is seen to be a $\CC$-linear endomorphism
$P \in \roman{End}_{\CC}{\CC}[X]$ such that
\linebreak
$(\roman{ad}\CC[X])^{m+1}P =0$.
 
Denote by $F_mD(X)$ the set of all linear differential operators on $X$
of order not greater than $m$. Clearly
$$ F_mD(X) \subset F_{m+1}D(X) \quad (m \geq 0)$$
and it is easily seen that $F_mD(X) \, F_lD(X) \subset F_{m+l}D(X)$.
Thus the set of all linear differential operators on $X$ forms a
$\CC$-algebra
$$ D(X) = \bigcup_{m=0}^{\infty} F_mD(X) $$
with filtration $F$. Note also that $F_0D(X) = {\Cal O}_X(X)$ by the
correspondence $P \mapsto P(1)$.
 
The sheaf $D_X$ of algebras of linear differential operators on $X$
is defined by the functor
$$
D_X : U \longmapsto D(U) \quad \text{for every open} \  U \subset X
$$
with obvious restriction maps. Thus $D_X(U) = D(U) =
\bigcup_{m=0}^{\infty} F_mD(U)$. The sheaf $D_X$ also has the increasing
filtration $F$ by orders $(F_mD_X)(U) \mathbreak
= F_mD(U) \  (m \geq 0)$.
 
The following lemma guarantees calculation in the algebraic case
similar to the complex analytic case.
 
\proclaim{Lemma}In a smooth n-dimensional algebraic variety $X$, every
point $p \in X$ has an affine open neighborhood $U$ with vector fields
$\partial_i$ and functions $x_i \  (1 \leq i \leq n)$ on $U$ satisfying
$$
[\partial_i, \partial_j] = 0, \quad [\partial_i, x_j] = \delta_{ij}
$$
$$
F_mD_X(U) = \bigoplus_{|\alpha| \leq m} {\Cal O}_X(U) \partial^{\alpha}
$$
where $\alpha = (\alpha_1, \cdots , \alpha_n)$ is a multi-index
$(|\alpha | = \sum_{i=1}^{n} \alpha_i)$ and
$$
\partial^{\alpha} = \prod_{i=1}^{n} \partial_{i}^{\alpha_i} .
$$
\endproclaim
 
\demo{Proof}Take $(x_i)$ to be a regular system at $p$, i.e., $\{ x_i \}$
generates the maximal ideal of ${\Cal O}_{X,p}$ and the differentials
$dx_i$ are linearly independent at $p$. There then exists an open $U$
such that
$$
f : U \longrightarrow \CC^n \quad (f(q) = (x_1(q), \cdots , x_n(q))
$$
is an etale map. The standard vector fields
$\displaystyle{\frac{\partial}{\partial z_i}}$ on $\CC^n$
lift uniquely to $\partial_i$ on $U$
$( df(\partial_i) = \displaystyle{\frac{\partial}{\partial z_i}})$
and $\{ x_i, \partial_i \}$ satisfies the requirement. In fact, for
$P \in F_mD_X(U)$ and $\alpha$ such that $|\alpha | = m$, put
$$
a_{\alpha}(x) = (-1)^m(\alpha !)^{-1}(\roman{ad}x_1)^{\alpha_1}
\cdots (\roman{ad}x_n)^{\alpha_n}P
$$
where $\alpha ! = \alpha_1 ! \cdots \alpha_n !$.
Then $a_{\alpha}(x)$ is of order 0 and hence $a_{\alpha}(x) \in {\Cal O}_X(U)$.
It is easily seen that $P - \sum_{|\alpha |=m} a_{\alpha}(x) \partial^{\alpha}$
is of order less than $m$. By induction, the lemma has been proved.
q.e.d.
\enddemo
 
\noindent
{\sl Remark.} Let $X_{\roman{an}}$ be the underlying complex manifold of a
smooth algebraic variety $X$ and $i:X_{\roman{an}} \rightarrow X$ the natural
morphism of local ringed spaces
($i^{-1}{\Cal O}_X \rightarrow {\Cal O}_{X_{\roman{an}}}$ is the identification
of regular functions on $X$ with holomorphic functions on $X_{\roman{an}}$).
Thus the sheaf $D_{X_{\roman{an}}}$ of linear differential operators with
holomorphic coefficients is regarded as
${\Cal O}_{X_{\roman{an}}} \bigotimes_{i^{-1}{\Cal O}_X} i^{-1}D_X$.
For a small open $U$ in $X_{\roman{an}}$ (in the classical topology)
the above choice of coordinates $\{ x_i, \partial_i \}$ is a standard one
in $D_{X_{\roman{an}}}(U)$.
 
\vskip 10mm

\noindent
{\bf 3. Filtrations of $D$-modules}
\bigskip
\noindent
3.1 Symbols.
 
\smallskip
The sheaf $D_X$ of algebras of linear differential operators has the
increasing filtration $F$ by orders, i.e., for an open $U$ in $X$,
$$
(F_mD_X)(U) = \{ P \in D_X(U) \  | \  \roman{ord}\, P \leq m \} .
$$
(Almost tautologically, $\roman{ord}\, P = m$ if and only if $P \in F_mD_X
\setminus F_{m-1}D_X$.)
Recall the following properties:
 
1) $F_mD_X \subset F_{m+1}D_X$,
 
2) $F_mD_XF_lD_X = F_{m+l}D_X$,
 
3) $F_mD_X$ is ${\Cal O}_X$-coherent,
 
4) $D_X = \displaystyle{\bigcup_{m=0}^{\infty}F_mD_X}$.
 
\noindent
Let $\roman{gr}\, D_X$ be the gradation of this algebra $D_X$
by the order filtrarion $F$,
$$
\roman{gr} \, D_X = \bigoplus_{m=0}^{\infty} \roman{gr}_m D_X
$$
where $\roman{gr}_m D_X = F_mD_X/F_{m-1}D_X$.
Note that for an affine open $U$,
$$
(\roman{gr}_m D_X)(U) = F_mD_X(U)/F_{m-1}D_X(U) .
$$
The graded algebra $\roman{gr} \, D_X$ is a commutative
${\Cal O}_X$-algebra since
$$
\roman{ord} \, [P, Q] \leq \roman{ord}(PQ)-1, \quad
(P, Q \in D_X).
$$
 
In the choice of local coordinates $\{x_i, \partial_i \}$ in Lemma in 2,
the projection
$$
F_mD_X(U) \longrightarrow \roman{gr}_m D_X(U)
$$
is realized as the symbol map
$$
P = \sum_{|\alpha | \leq m} a_{\alpha} \partial^{\alpha}
\longmapsto
\sum_{|\alpha |=m} a_{\alpha} \xi^{\alpha} = \sigma_m(P)
$$
where $\xi = (\xi_1, \cdots , \xi_n)$ is the linear coordinate
system corresponding to $(x_i)$ on the cotangent bundle
$T^*U$. Thus the symbol $\sigma_m(P)$ is regarded as an element of
the polynomial algebra ${\Cal O}_X(U)[\xi_1, \cdots , \xi_n]$ over
${\Cal O}_X(U)$. It is a standard fact that this symbol map
$\sigma$ is independent of the choice of coordinates and it gives
rise to the following global identification of the graded algebras:
$$
\roman{gr} \, D_X \, \widetilde{\longrightarrow} \, \pi_*{\Cal O}_{T^*X}
$$
where $\pi : T^*X \rightarrow X$ is the cotangent bundle of $X$
and $\pi_*$ is the operation of a direct image sheaf
($(\pi_*{\Cal O}_{T^*X})(U) = \Gamma (\pi^{-1}(U), {\Cal O}_{T^*X})
= {\Cal O}_X(U)[\xi_1, \cdots , \xi_n]$).
 
\bigskip
\noindent
3.2 Good filtrations.
 
\smallskip
A left $D_X$-module $M$ simply means a sheaf of left $D_X$-modules:
a sheaf $M$ on $X$ such that for every open $U$ in $X$, $M(U)$ is a left
$D_X(U)$-module compatible with restriction data. A right $D_X$-module
is similarly defined. Since ${\Cal O}_X$ is a subalgebra of $D_X$,
a $D_X$-module has the natural structure as an ${\Cal O}_X$-module.
On an algebraic variety $X$, we usually consider $D_X$-modules
which are ${\Cal O}_X$-{\it quasi-coherent} in order to pursue smooth
manipulation in algebraic geometry. (An ${\Cal O}_X$-module $F$ is
called quasi-coherent if $F|_U$ is isomorphic to the sheaf made by
localization of the ${\Cal O}_X$-module $F(U)$ on every affine open $U$.)
 
However, we retain that for deeper analysis of solutions to equations,
one often needs non-quasi-coherent ${\Cal O}_X$-modules. In particular,
for analysis on the complex manifolds $X_{\roman{an}}$, these are
sometimes essential, e.g., a sheaf of distributions and/or
hyperfunctions etc,...
 
At any rate, usual systems of linear partial differential equations
correspond to more restricted $D_X$-modules, $D_X$-{\it coherent}
modules. Here we take the definition of the $D_X$-coherency as that of
the local finite presentation of $D_X$-modules, i.e., we call a $D_X$-module
$M$ on $X$ $D_X$-coherent if every point of $X$ has a neighborhood
$U$ with an exact sequence
$$
D_U^m \longrightarrow D_U^l \longrightarrow M|_U \longrightarrow 0
$$
where $D_U$ is the sheaf of algebras of linear differential operators
on $U$ considered as a $D_U$-module.
 
As is earlier introduced, the category of coherent $D_X$-modules
corresponds to that of systems of linear partial differential
equations on $X$. Note that on an affine open $U$, $M|_U$ is the
localization of $M(U)$ and hence the above condition corresponds to the
exact sequence of $D(U)$-modules:
$$
D(U)^m \overset\Phi\to\longrightarrow D(U)^l 
\overset\Psi\to\longrightarrow M(U) \longrightarrow 0
$$
 
\smallskip
\noindent
{\sl Remark.} The above coherent $D$-module corresponds to the system of
linear partial differential equations
$$
\sum_{j=1}^{l} P_{ij} \, u_j = 0 \quad (1 \leq i \leq m),
$$
by the maps
$$
\Phi (Q_1, \cdots , Q_m) = (Q_1, \cdots , Q_m)(P_{ij}),
$$
$$
\Psi (R_1, \cdots , R_l) = \sum_{j=1}^{l} R_j \, u_j
$$
where $(P_{ij})$ is an $m \times l \,$-matrix of entries in $D(U)$.
 
In order to define the characteristic variety of a $D$-module,
we introduce the notion of good filtrations matching the order
filtration of $D_X$.
 
For a coherent $D_X$-module $M$, let $F_mM \subset M (m \in {\Bbb Z})$
be an increasing filtration by coherent ${\Cal O}_X$-submodules (i.e.,
${\Cal O}_X$-submodules of locally finite presentation) such that
 
1) $F_mM = 0 \quad (m \ll 0)$,
 
2) $\displaystyle{M = \bigcup_{m \in {\Bbb Z}} F_mM}$,
 
3) $F_lD_X \, F_mM \subset F_{l+m}M \quad (l, m \in {\Bbb Z})$.
 
\noindent
Then by 1), 2), 3), the gradation of $M$ by $F$
$$
\roman{gr}^F M = \bigoplus_{m \in {\Bbb Z}} F_mM/F_{m-1}M
$$
is a graded $\roman{gr}\, D_X$-module.
 
\smallskip
\noindent
{\bf Definition.} A filtration $F$ of a coherent $D_X$-module $M$ is called
{\it good} if
$$
F_lD_X \, F_mM = F_{l+m}M  \quad \text{for } m
\text{ large enough and all } l \geq 0.
$$
Here the left hand side is the ${\Cal O}_X$-submodule generated by the
multiplication of $F_mM$ by $F_lD_X$.
 
\smallskip
The following is then rather easily proved ([11, II, Prop.1.2.3]).
 
\proclaim{Proposition}
A filtration $F$ is good if and only if the graded module
$\roman{gr}^FM$ is a coherent
$\roman{gr} \, D_X \, (= \pi_*{\Cal O}_{T^*X})$-algebra.
\endproclaim
 
By definition, a coherent $D_X$-module $M$ is locally finitely generated, i.e.,
on an open $U$
$$
M|_U = \sum_{i=1}^{l}D_U \, u_i \quad (u_i \in M(U)).
$$
Define $F_mM|_U = \sum_{i=1}^{l} F_mD_U \, u_i \  (m \geq 0)$. Then $F_m$ is
a good filtration of $M|_U$. Thus any coherent $D_X$-module {\it locally}
has a good filtration. In the algebraic case, any coherent $D_X$-module
{\it globally} has a good filtration thanks to the quasi-compactness of the
Zariski topology (see [11, II,1.2]).
 
\smallskip
\noindent
{\bf Example.} Let $M = D_X \, u$ ($D_X$-cyclic by a section $u \in M(X)$).
$F_mM = F_mD_X \, u$ then gives a good filtration of $M$. Put
$$
I = \roman{Ann} \, u = \{ P \in D_X \, | \, Pu = 0 \, \}
$$
the annihilator of $u$ (thus $M \simeq D_X/I$). Then
$$
\roman{gr}^FM \simeq \roman{gr}\, D_X/ \roman{gr}\, I
$$
where
$$
\roman{gr}\, I = \bigoplus_{m \leq 0} F_mI/F_{m-1}I
\subset \pi_*{\Cal O}_{T^*X},
\quad (F_mI = I \cap F_mD_X).
$$
 
\vskip 10mm
 
\noindent
{\bf 4. Characteristic varieties}
 
\bigskip
\noindent
4.1 Definition.
 
\smallskip
Let $M$ be a coherent $D_X$-module and $F$ its good filtration. Then the
graded $\pi_*{\Cal O}_{T^*X}$-module $\roman{gr}^FM$ is coherent where
$\pi : T^*X \rightarrow X$ is the cotangent bundle. Since the sheaf pull-back
$\pi^{-1} \roman{gr}^FM$ on $T^*X$ is $\pi^{-1} \pi_*{\Cal O}_{T^*X}$-coherent,
we have an ${\Cal O}_{T^*X}$-coherent module
$$
\pi^{\bullet} (\roman{gr}^FM) = {\Cal O}_{T^*X} \otimes_
{\pi^{-1}\pi_*{\Cal O}_{T^*X}}\pi^{-1}(\roman{gr}^FM)
$$
on the cotangent bundle $T^*X$ (the algebra homomorphism
$\pi^{-1}\pi_*{\Cal O}_{T^*X} \mathbreak
\rightarrow {\Cal O}_{T^*X}$ is a natural
restriction of functions). The {\it characteristic variety}
$\roman{ch} \, M$ of $M$ is then defined to be the support of
the ${\Cal O}_{T^*X}$-coherent module
$\pi^{\bullet}(\roman{gr}^FM)$. By the coherency, the characteristic
variety is an algebraic subvariety of $T^*X$ conic along the fibers
(= cotangent spaces).
 
We shall look at it in a more naive way. Let $U$ be a small affine
open set in $X$. Then
$$
(\pi_*{\Cal O}_{T^*X})(U) = {\Cal O}_X(U) \otimes_{\CC}
\CC [\xi_1, \cdots , \xi_n]
$$
where $(\xi_1, \cdots , \xi_n)$ is a coordinate system of the cotangent
space (we assume $T^*U \simeq U \times \CC^n$). For $\roman{gr}^FM$,
we have
$$
(\roman{gr}^FM)(U) = \bigoplus_{m \in {\Bbb Z}} F_mM(U)/F_{m-1}M(U),
$$
which is a graded module over the graded algebra
${\Cal O}_X(U) \otimes_{\CC} \mathbreak
\CC [\xi_1, \cdots , \xi_n]$.
The ${\Cal O}_{T^*X}$-module $\pi^{\bullet}(\roman{gr}^FM)$ is
simply the localization of the
${\Cal O}_X(U) \mathbreak \otimes_{\CC} \CC [\xi_1, \cdots , \xi_n]$-module
$(\roman{gr}^FM)(U)$ on $T^*U \simeq U \times \CC^n$ and hence the
characteristic variety $\roman{ch}\, M$ on $T^*U$ is nothing but the zeroes
(affine subvariety) of the annihilator ideal in the algebra
${\Cal O}_X(U) \otimes_{\CC} \CC [\xi_1, \cdots , \xi_n]$
$$
\roman{Ann}_{{\Cal O}_X(U) \otimes_{\CC} \CC [\xi_1, \cdots , \xi_n]}
(\roman{gr}^FM)(U) .
$$
By the gradedness, this variety is conic along the fibers $\CC^n$.
 
\smallskip
\noindent
{\bf Example.} Let $M = D_X \, u$ be as in Example in 3.2. Then
$\roman{ch}\, M = V(\roman{gr}\, I)$ the zeroes defined by the ideal
$\roman{gr}\, I \subset \pi_*{\Cal O}_{T^*X}$. Notice that even if
$P_1, \cdots , P_m$ are generators of the ideal $I$ in $D_X$,
the symbols $\sigma(P_1), \cdots , \sigma(P_m)$ are not necessarily
generators of $\roman{gr}\, I$. That is, the characteristic variety
$\roman{ch}\, M$ is not exactly the zeroes of the symbols
$\sigma(P_i) \  (1 \leq i \leq m)$ (is contained in those).
Here we see the importance of the concept of $D_X$-modules in defining
the characteristic varieties of systems of linear partial differential
equations. It is however known that there exist generators
$P_i$'s in $I$ such that the symbols $\sigma(P_i)$'s generate
$\roman{gr}\, I$ (see [11, II,2]).
 
\proclaim{Theorem} The characteristic variety $\roman{ch} \, M$ is
independent of the choice of good filtrations of a coherent
$D_X$-module $M$.
\endproclaim
 
For the proof, see [11, II,Th.2.1].
By this theorem, the characteristic variety turns out to be a true
invariant of a $D$-module. Also in the analytic case, since a
coherent $D$-module locally has a good filtration, the characteristic
variety can be defined globally by a similar theorem.
 
\bigskip
 
\noindent
4.2 The fundamental theorem.
 
\smallskip
 
As is well-known, the cotangent bundle $T^*X$ has the canonical
symplectic structure $\omega$ which is expressed in local coordinates
$$
\omega = \sum_{i=1}^{n} d\xi_i \wedge dx_i \quad (n = \dim X).
$$
The symplectic structure $\omega$ defines the Poisson bracket in
the space of functions on $T^*X$ which is expressed in local coordinates
$$
\{ f, g \} = \sum_{i=1}^{n} (\frac{\partial f}{\partial \xi_i}
\frac{\partial g}{\partial x_i} - \frac{\partial g}{\partial \xi_i}
\frac{\partial f}{\partial x_i})
\quad (f, g \in {\Cal O}_{T^*X}).
$$
 
A subvariety $V$ in $T^*X$ is said to be {\it involutive} if the defining
ideal $I(V) \in {\Cal O}_{T^*X}$ is closed under the Poisson Bracket
$\{\  , \  \}$. It is easily seen that if $V$ is involutive, then the
tangent space of a smooth point $p \in V$ is involutive with respect to
the symplectic form $\omega$ ($T_pV^{\bot} \subset T_pV$ where $\bot$
denotes the orthogonal complement with respect to
$\omega$ in $T_p(T^*X)$) and hence every irreducible
component of $V$ has dimension not less than
$\dim X = \displaystyle{\frac{1}{2}} \dim T^*X \  (\text{if } V
\not= \emptyset)$.
 
The following theorem is called the fundamental theorem of
algebraic analysis, which is first proved in [23]. Later O. Gabber gave
a purely algebraic proof for this theorem [5] but still difficult.
 
\proclaim{Theorem (Sato-Kawai-Kashiwara)}The characteristic variety
of a coherent $D$-module is involutive.
\endproclaim
 
\smallskip
In particular, if $M \not= 0 \  (M = 0 \Leftrightarrow \roman{ch}\, M =
\emptyset)$, then
$$
\dim (\text{component of ch}\, M) \geq \dim X.
$$
Note that the weaker statement $`` \dim \, \roman{ch}\, M \geq \dim X"$
is much more easily proved (in the algebraic case) by J. Bernstein
(see [1], [11, II,Th.5.1]).
 
Now we shall ask the following question. What kinds of non-zero coherent
$D$-modules are of ``smallest size"? In a monogenic case $M = D_Xu$,
those must correspond to the case when the annihilator ideals
$\roman{Ann}_{D_X}u$ have
``largest size" which means that the corresponding systems of linear
partial differential equations are ``maximally overdetermined". Taking
the gradation $\text{gr Ann}\, u$, the characteristic variety
$\roman{ch}\, M = V(\text{gr Ann}\,u)$ must have ``smallest size" as
possible but these are of dimension not less than $\dim X$.
Thus we attain the case of holonomic $D$-modules.
 
\smallskip
\noindent
{\bf Definition.} A coherent $D_X$-module $M$ is called {\it holonomic}
if $\dim \roman{ch}\, M \break = \dim X$ or $M = 0$.
 
\smallskip
All substantial examples in these notes are holonomic.
 
\vskip 10mm
 
\noindent
 
\noindent
{\bf 5. Examples}
 
\bigskip
\noindent
5.1. Ordinary differential equations.
\smallskip
 
Let $P(x, \partial ) = \sum_{i=0}^{m} a_i(x) \partial^i$ be a
non-zero linear operator on $\CC \  (a_i(x) \break \in \CC [x], \partial
= \displaystyle{\frac{d}{dx}}, a_m(x) \not= 0)$.
The ordinary differential equation
$$
P(x, \partial ) \, u = 0
$$
corresponds to the $D_{\CC}$-module
$$
M = D_{\CC} \, u = D_{\CC}/I \quad (I = D_{\CC} \, P(x, \partial )) .
$$
In this case, $\roman{gr} \, I$ turns out to be a principal ideal
generated by the symbol $\sigma_m(P)(x, \xi ) \mathbreak = a_m(x) \xi^m
\in \CC [x, \xi ]$ and hence
$$
\align
\roman{ch} \, M &= \{ (x, \xi ) \in \CC \times \CC \, | \,
a_m(x) \xi^m = 0 \, \} \\
 &= \CC \times 0 \cup (\bigcup_{a_m(x)=0} x \times \CC)
\endalign
$$
i.e., the union of the zero section and the fibers at the singular
points $a_m(x) = 0$ ($T^* \CC = \CC \times \CC$). Hence
$\dim \roman{ch}\, M = 1$ and $M$ is holonomic.
 
\bigskip
\noindent
5.2 Connections.
 
A $D_X$-module $M$ is called a {\it connection} if $M$ is a locally
free ${\Cal O}_X$-module of finite rank (i.e., vector bundle as an
${\Cal O}_X$-module).
 
Let $M$ be a connection and $U$ open in $X$. Let $\Theta_X$ be the sheaf
of vector fields on $X$. Write the action of a vector field $\theta
\in \Theta_X(U)$ on a section $s \in M(U)$ as
$$
\nabla_{\theta} s = \theta s.
$$
Then by the definition of the $D_X(U)$-action, we have
 
i) $\nabla_{f \theta}s = f \nabla_{\theta}s \quad (f \in {\Cal O}_X(U))$,
 
ii) $\nabla_{\theta}(fs) = \theta (f) s + f \nabla_{\theta} s
\quad ([\theta , f] = \theta (f) )$,
 
iii) $\nabla_{[\theta, \theta ']} s = [\theta , \theta '] \, s$.
 
\noindent
The conditions i), ii) are the usual ones for a ``connection"
and iii) corresponds to the ``integrability". Thus a connection in our sense
is an integrable connection.
 
Conversely when a vector bundle $M$ is given with $\Theta_X$-action
through $\nabla$ satisfying i), ii), iii), then this action extends
to the left $D_X$-action and $M$ acquires the structure of a $D_X$-module.
 
\smallskip
\noindent
{\bf Example.} Let $M = D_{\CC} \, u = D_{\CC}/D_{\CC}P $ be as in 5.1.
Then $M|_U$ is a connection on the open set
$\{ x \in \CC \, | \, a_m(x) \not= 0 \, \}$. In fact, let
$u_i = \partial^iu \  (0 \leq i < m)$. Then
$$
M|_U = D_Uu \simeq \sum_{i=0}^{m-1} {\Cal O}_Uu_i
$$
since $a_m(x)^{-1} \in {\Cal O}(U)$.
 
\medskip
Let $M$ be a connection on $X$. Define a filtration $F$ on $M$ by
$$
F_iM = 0 \quad (i \leq 0), \qquad F_iM =M \quad (i \geq 1).
$$
Since $M$ is ${\Cal O}_X$-coherent, this is clearly a good filtration.
Since
$$
\roman{gr}^FM = M \quad \text{(concentrated at degree 1)},
$$
$\roman{gr}_1D_X$ then acts as 0 on $\roman{gr}^F_1 M = M$.
This means, for instance, that locally every linear coordinate
$\xi_i$ on $T^*X$ belongs to $\roman{Ann} \, \roman{gr}^FM$,
which implies
$$
\roman{ch} \, M = T^*_XX \  (= \text{zero-section of }
\pi : T^*X \rightarrow X).
$$
 
Actually we know:
 
\proclaim{Theorem} For a coherent $D_X$-module $M$,
the following are equivalent.
 
1) $M$ is a connection.
 
2) $M$ is ${\Cal O}_X$-coherent.
 
3) $\roman{ch} \, M = T^*_XX.$
 
\endproclaim
 
For the proof, see [11, II,Prop.2.3].
 
We close this section by citing an important classical theorem
in the {\it analytic} case.
 
\proclaim{Theorem (Frobenius)} Let $M$ be an {\it analytic} connection
on a {\it complex} manifold $X$. Define the subsheaf of vector spaces
in $M$:
$$
(DR \, M)(U) = \{ s \in M(U) \, | \nabla_\theta s = 0
\text{ for any vector field } \theta \in \Theta_X(U) \, \}
$$
for every open $U$ in $X$. Then $DR \, M$ is a local system of
rank equal to $\roman{rank} \, M$ (i.e., for a small connected $U$,
$(DR \, M)(U) \simeq \CC^r$ where $r = \roman{rank} \, M \  
(M|_U \simeq {\Cal O}_U^r)$).
 
\endproclaim
 
On a complex manifold $X$, this correspondence
 
\centerline{
\{ analytic connections on $X$ \} $\overset DR\to\longrightarrow$
\{ local systems on $X$ \}
}
 
\noindent
gives rise to a categorical equivalence. The quasi-inverse of the
functor $DR$ (the de Rham functor) is given by
$$
L \longmapsto {\Cal O}_X \otimes_{\CC} L
$$
where the connection in the right-hand-side is defined by
$\nabla_\theta(f \otimes s) = \theta (f) \otimes s
\quad (\theta \in \Theta_X, f \in {\Cal O}_X, s \in L)$.
 
In the algebraic case, the above correspondence $DR$ causes much delicate
problems. Now let $X$ be a smooth algebraic variety and $X_{\roman{an}}$
the underlying complex manifold. For an algebraic connection $M$ on $X$,
in turn, define the functor $DR$ by
$$
DR \, M = DR \, M_{\roman{an}}
\quad \text{where } M_{\roman{an}} =
{\Cal O}_{X_{\roman{an}}} \otimes_{{\Cal O}_X}M .
$$
Thus $DR$ is a functor from the category of algebraic connections on $X$
into that of local systems on the complex manifold $X_{\roman{an}}$.
In this case, $DR$ does not give rise to a categorical equivalence but
by restricting the class of algebraic connections, P. Deligne established
a nice equivalence between these categories ([4]). For this, the
notion of regularities is necessarily involved, and an algebraic
connection is called {\it regular} if its (unique) meromorphic
extension to a compactification $\overline X$ has only regular
singularities at the boundary $\overline X \setminus X$.
(For the precise definition, see [11, IV].)
 
\proclaim{Theorem (Deligne)} The functor $DR$ gives a categorical
equivalence between the category of regular connections on an algebraic
variety X and that of local systems on the complex manifold
$X_{\roman{an}}$.
 
\endproclaim
 
This correspondence $DR$ is intensively generalized to $D$-modules
and plays a substantial role in the Riemann-Hilbert correspondence
for regular holonomic $D$-modules (M. Kashiwara and Z. Mebkhout,
see [16], [20], [11, V]).
 
\noindent
{\bf Example.} For a fixed $\lambda \in \CC$, let ${\Cal O}_\lambda$
be the Euler system
$$
{\Cal O}_\lambda = D_{\CC} u = D_{\CC} / D_{\CC} (x \partial - \lambda )
\quad \text{ on } \CC \  (\partial = \frac{d}{dx}).
$$
Then ${\Cal O}_\lambda |_{\CC^{\times}} \  (\CC^{\times} =
\CC \setminus \{ 0 \})$ is an algebraic connection and the
anaytic solution sheaf $S_\lambda$ is a local system of rank one
generated by the multi-valued holomorphic function $x^\lambda$ on
$\CC^\times$
 
On the other hand,
$$
DR({\Cal O}_\lambda |_{\CC^{\times}}) =
\CC x^{- \lambda} u \subset {\Cal O}_\lambda |_{\CC^{\times}}.
$$
since
$$
x \partial (x^{- \lambda}u) = - \lambda x^{- \lambda}u + x^{- \lambda}
(x \partial u) = 0
$$
by $x \partial u = \lambda u$. In this sense, $S_\lambda$ and
$DR({\Cal O}_\lambda |_{\CC^{\times}})$ are local systems dual to
each other. Further we see, even for an algebraic $D_{\CC}$-module
${\Cal O}_\lambda$ as above, this correspondence makes sense only in the
analytic category (we need such analytic functions like $x^\lambda$).
 
\vskip 10mm
 
\noindent
{\bf 6. Basic functors}
 
\bigskip
 
\noindent
6.1 Inverse images.
 
\smallskip
Let $X \to Y$ be a morphism of smooth algebraic varieties. For a
$D_Y$-module $M$, let
$$
f^*M = {\Cal O}_X \otimes_{f^{-1}{\Cal O}_Y}f^{-1}M
$$
be the inverse image as an ${\Cal O}$-module. Then $f^*M$ turns out to
be a $D_X$-module by the chain rule. That is, if $\{y_i, \partial_i \}$
is a local coordinate system on $Y$, vector fields $\theta \in \Theta_X$
act on $f^*M$ by
$$
\theta (\phi \otimes u) = \theta (\phi ) \otimes u +
\phi \sum_{i=1}^{n} \theta (y_i \circ f) \otimes \partial_i u
\quad (\phi \in {\Cal O}_X, u \in M)
$$
and these actions extend to the $D_X$-action.
 
Even if $M$ is $D_Y$-coherent, $f^*M$ is not necessarily $D_X$-coherent.
For example, on $Y = \CC^2$ consider an equation
$$
\partial_yu = 0 \quad (y \text{ is the second coordinate}).
$$
The corresponding $D_Y$-module $M \simeq D_Y/D_Y \partial_y \simeq
D_{\CC} \boxtimes {\Cal O}_{\CC}$ has the inverse image
$$
i^*M \simeq (\sum_{j=0}^{\infty} \CC \partial_x^j) \otimes_{\CC}
{\Cal O}_{\CC}
$$
for an inclusion $i : \CC \hookrightarrow \CC^2 \  (i(y) = (0, y))$.
In the above $\boxtimes$ denotes the outer tensor product (over $\CC$)
on the product variety. Thus, as an $D_{\CC}$-module, $i^*M \simeq
\displaystyle{\bigoplus_{j=0}^{\infty}{\Cal O}_{\CC}}$
is not a finitely generated $D_{\CC}$-module.
 
However, it is known that if $M$ is holonomic then so is $f^*M$
(see [11, III]).
 
\bigskip
\noindent
6.2 Direct images (integrations along fibers).
\smallskip
 
In contrast with the inverse images, the definition of the ``direct
image" of $D$-modules is rather complicated. We want a certain
$D_Y$-module ``$f_\bullet M$" on $Y$ for $f:X \to Y$ and for a $D_X$-module
$M$.
 
As an extreme example, consider a closed immersion
$$
i : X \longrightarrow X \times \CC = Y \quad (i(x) = (x, 0)).
$$
For a $D_X$-module $M$ on $X$, let $i_*M$ be the direct image
in the sheaf theory:
$$
(i_*M)(U) = M(U \cap (X \times 0)).
$$
This becomes an ${\Cal O}_Y$-module as usual but not a $D_Y$-module.
Since $D_Y \simeq D_X \boxtimes D_{\CC} \simeq
\displaystyle{\bigoplus_{j=0}^{\infty} D_X \boxtimes {\Cal O}_{\CC}
\partial^j} \  (\partial = \displaystyle{\frac{\partial}{\partial t}},
\, t \text{ :coordinate of } \CC)$, we want a $\partial$-action on some
$i_\bullet M$, an extension of $i_*M$. Considering the infinite sum
$$
\bigoplus_{j=0}^{\infty} i_*M \otimes \partial^j
\quad (i_*M \otimes \partial^j \simeq i_*M 
\text{ as } D_X \text{-modules})
$$
with the $\partial$-action by $\partial (u \otimes \partial^j) =
u \otimes \partial^{j+1}$ (${\Cal O}_{\CC}$ -action by
$\phi (t) u \otimes \partial^j = \phi (0) u \otimes \partial^j$),
$i_*M$ extends to a $D_{X \times \CC}$-module.
 
In general, we have to be much more careful. So far we have only
considered {\it left} $D$-modules since they correspond naturally
to systems of linear partial differential equations. But in order to
understand the ``direct images", it is more convenient to consider
{\it right} $D$-modules. These simply correspond to the adjoint systems
of usual ones.
 
For example, let $P^*$ be the adjoint operator of $P \in D_{\CC^n}
\  (\partial_i \mapsto - \partial_i)$. Then $(PQ)^* = Q^*P^*$. If $M$
is a left $D_{\CC^n}$-module, by the action
$$
u \, P = P^* u \quad (u \in M, P \in D_{\CC^n})
$$
$M$ is regarded as a {\it right} $D_{\CC^n}$-module.
 
Globally in general, this procedure can be defined by using Lie
derivatives on highest differential forms. Let $\Omega_X$ be the sheaf
of highest differential forms on $X$ (the canonical line bundle).
The {\it Lie derivative} $L_\theta \omega$ for a vector field
$\theta \in \Theta_X$ and $\omega \in \Omega_X$ is by definition
$$
(L_\theta \omega)(\theta_1, \cdots , \theta_n) =
\theta (\omega (\theta_1, \cdots , \theta_n)) -
\sum_{i=1}^{n} \omega (\theta_1, \cdots , [\theta , \theta_i], \cdots,
\theta_n)
$$
where $\theta_i \in \Theta_X$.
Then $L_{\phi \theta} \omega = L_\theta (\phi \omega)$ for
$\phi \in {\Cal O}_X$. Hence defining the $\Theta_X$-action on
$\Omega_X$ by $\omega \theta = - L_\theta \omega$, $\Omega_X$
gains the right $D_X$-module structure.
 
For a {\it left} $D_X$-module $M$, $M \otimes_{{\Cal O}_X} \Omega_X$
then turns out to be a {\it right} $D_X$-module by
$$
(u \otimes \omega ) \theta = - (\theta u) \otimes \omega +
u \otimes \omega \theta \quad (u \otimes \omega \in M \otimes \Omega_X,
\theta \in \Theta_X).
$$
 
In the above example of $D_{\CC^n}$-modules, we have fixed a global
section $dx_1 \wedge \cdots \wedge dx_n$ of $\Omega_{\CC^n}$ and
identify $M$ with $M \otimes_{{\Cal O}_{\CC^n}} \Omega_{\CC^n}$.
 
Now let $f : X \to Y$ be a morphism of smooth algebraic varieties.
As is seen in 6.1, the inverse image $f^*D_Y$ (as an ${\Cal O}$-module)
ia a left $D_X$-module. Simultaneously it is also a right
$f^{-1}D_Y$-module commuting with the left $D_X$-action. We write
$$
D_{X \to Y} = f^*D_Y
$$
as a double $(D_X , f^{-1}D_Y)$-module.
 
If $M$ is a {\it right} $D_X$-module on $X$, then
$M \otimes_{D_X} D_{X \to Y}$ turns out to be a right $f^{-1}D_Y$-module
(coming from the right action on $D_{X \to Y}$).
Hence if we take the sheaf direct image
$$
f_*(M \otimes_{D_X} D_{X \to Y}),
$$
then it becomes a {\it right} $D_Y$-module on $Y$.
 
For a {\it left} $D_X$-module $M$, we apply the left-right correspondence
(tensoring $\otimes \Omega_X$) and take the above procedure. Thus the
final form becomes the following messy one
$$
f_*((M \otimes_{{\Cal O}_X} \Omega_X ) \otimes_{D_X} D_{X \to Y}
\otimes_{f^{-1}{\Cal O}_Y} f^{-1}(\Omega _Y^{-1}))
$$
a left $D_Y$-module (tensoring $\Omega_Y^{-1}$ converts the right
$D_Y$-structure into the left one). Furthermore, the description
of the left $D_Y$-action is highly complicated in general.
 
\smallskip
\noindent
{\sl Exercise.} Check the earlier example for $X \hookrightarrow
X \times \CC$.
 
\smallskip
But still the above definition is not a final one in general
(if $X, Y$ are not affine). For the correct definition of the direct
images of $D$-modules, we have to use a concept of derived categories
and several operations in them. We are not going into details now,
but only write down
$$
f_\bullet M = \int_f M = Rf_*((M \otimes_{{\Cal O}_X} \Omega_X)
\overset L\to\otimes_{D_X} (D_{X \to Y}
\otimes_{f^{-1}{\Cal O}_Y} f^{-1}(\Omega_Y)^{-1}))
$$
(see [11, I]).
 
This functor works well in the ``Riemann-Hilbert" correspondence
for regular holonomic $D$-modules and corresponds simply to
the sheaf direct image functor for the solution sheaves or the
de Rham complexes.

 
 
 
 
\newpage
 
\centerline{\bf II. Equivariant $D$-modules, examples}
 
\vskip 10mm
\noindent
{\bf 1. Definition}
 
\bigskip
 
When a variety has a group action, it is often important to consider
a class of systems of linear partial differential equations with
equivariance property under the group action. Here we start with defining
the group equivariance property matching the $D$-module operations.
 
Let $G$ be an algebraic group acting on a smooth algebraic variety
$X$ and $\alpha : G \times X \to X$ the group action morphism
$(\alpha (g,x) = gx \  (g \in G, x \in X))$. Naively, for instance, a
sheaf $F$ on $X$ is considered as $G$-equivariant if it is given a datum
of sheaf morphisms $F_x \simeq F_{gx}$ ``simultaneously" for
$(g, x) \in G \times X$. More strictly, $F$ is said to be $G$-equivariant
if there exists a morphism
$$
\alpha^{-1}F \  \widetilde{\longrightarrow} \  \roman{pr}_X^{-1}F
$$
satisfying the associativity condition coming from the group
multiplication of $G$ (cocycle condition). ($\alpha^{-1},
\  \roman{pr}_X^{-1}$ are the inverse image functor in the sheaf theory
and $\roman{pr}_X : G \times X \to X$ is the projection onto $X$.)
 
We follow an analogous approach in $D$-module operations but extend the
setup a little wider, i.e., attach a twisting datum on the group $G$.
This extension contains much more examples, in particular, Sato's
relative invariants on prehomogeneous vector spaces and Gelfand's
generalized hypergeometric equations.
 
A twisted datum is a connection $L$ on $G$. We want to define a
$D_X$-module $M$ to be {\it `` L-twistedly" $G$-equivariant} if there
exists a $D$-module homomorphism (with naturality condition)
$$
\alpha^*M \simeq L \boxtimes M \quad \text{on } G \times X.
$$
Here $\alpha^* M$ is the inverse image of the $D$-module $M$ and
$L \boxtimes M$ is the outer tensor product over $\CC$ on the product
space $G \times X \  (L \boxtimes M = \roman{pr}_G^{-1}L \otimes_{\CC}
\roman{pr}_X^{-1}M \simeq \roman{pr}_G^*L \otimes_{{\OO}_{G \times X}}
\roman{pr}^*_XM)$.
 
But then in the diagram
$$
\CD
G \times G \times X @>{\mu \times 1_X}>> G \times X \\
@V{1_G \times \alpha}VV   @VV{\alpha}V \\
G \times X @>\alpha>> X
\endCD
$$
($\mu : G \times G \to G$ is the group multiplication),
the following are required by the naturality condition
$$
\mu^*L \boxtimes M \simeq L \boxtimes \alpha^*M.
$$
The right hand side is isomorphic to
$L \boxtimes L \boxtimes M$ by the above equivariance.
Hence our connection $L$ on $G$ should satisfy
$$
\mu^* L \simeq L \boxtimes L,
$$
i.e., $L$ itself should be $L$-twistedly $G$-equivariant on $G$
under the group multiplication $\mu$.
 
Furthermore in this case $L$ should necessarily be of rank one
(line bundle).
In fact, consider the maps
$$
G \times e \overset\iota \to\hookrightarrow G \times G
\overset\mu \to\longrightarrow G \quad
(e \in G \text{ : identity}).
$$
Then $\mu^* L \simeq L \boxtimes L$. Applying $\iota^*$, we have
$$
L \simeq \iota^* \mu^* L
$$
by naturality. But then the right hand side
$\simeq L \boxtimes e^* L$ where $e^*L$ is the geometric stalk of $L$
at $e$. Hence $e^*L \simeq \CC$, i.e., $L$ is a line bundle over $G$.
 
\smallskip
\noindent
{\bf Example.} Let $\frak g = \roman{Lie} \, G$
be the Lie algebra of $G$. Fix a Lie algebra homomorphism
$\lambda : \frak g \to \CC$ where $\CC$ is considered as an abelian
Lie algebra. For $\theta \in \frak g$, let $L_\theta$ be the
corresponding right invariant vector field on $G$, i.e.,
$(L_\theta f)(x) = \displaystyle{\frac{d}{dt}}f(e^{-t \theta} x)|_{t=0}$.
The system of linear partial differential equations
$$
(L_\theta - \lambda (\theta ))u = 0 \quad (\theta \in \frak g)
$$
corresponds to the $D_G$-module
$$
\OO_\lambda = D_G / \sum_{\theta \in \frak g} D_G (L_\theta - \lambda
(\theta)).
$$
Then $\OO_\lambda \simeq \OO_G u$ is a rank one connection with
$$
\nabla_\theta (fu) = \theta (f)u + \lambda (\theta )fu
\quad (\theta \in \frak g )
$$
and $\OO_\lambda$ is $\OO_\lambda$-twistedly $G$-equivariant.
 
\medskip
\noindent
{\sl Examples of Example.}
\smallskip
1) Let $G$ be a linear group and $\det : G \to \CC^\times$.
A multi-valued function $\det(x)^\lambda \  (x \in G)$ generates a rank
one connection
$$
D_G \det(x)^\lambda = \OO_G \det (x)^\lambda.
$$
This is a simplest connection with regular singularities.
 
\smallskip
2) On the additive group $\Bbb G_a = \CC$, a similar connection
$$
D_{\CC}e^{\lambda x} = D_{\CC}/D_{\CC}(\displaystyle{\frac{d}{dx}} -
\lambda )
$$
has an irregular singularity at $x = \infty$.
 
\medskip
Finally, we arrive at the definition of equivariant $D$-modules.
 
\smallskip
\noindent
{\bf Definition.} Let $X$ be a smooth algebraic variety with
algebraic group action $\alpha : G \times X \to X$. Let $L$ be a rank
one connection on $G$. A $D_X$-module $M$ is called $L$-{\it twistedly
$G$-equivariant} if there exists a $D$-module isomorphism
$$
\phi : \alpha^*M \  \widetilde{\longrightarrow} \ 
L \boxtimes M \quad \text{on } G \times X
$$
satisfying the following cocycle condition:
$$
\roman{pr}_{G \times X}^*(\phi ) \circ (1_G \times \alpha )^*(\phi )
\simeq (\mu \times 1_X)^*(\phi )
$$
with respect to the commutative diagram
$$
\CD
G \times G \times X @>{\mu \times 1_X}>> G \times X \\
@V{1_G \times \alpha}VV   @VV{\alpha}V \\
G \times X @>\alpha>> X.
\endCD
$$
Here $L$ itself necessarily turns out to be an $L$-twistedly
$G$-equivariant $D_G$-module under the group multiplication
$\mu : G \times G \to G$.
 
When $L = \OO_\lambda$ for a Lie algebra character
$\lambda : \frak g \to \CC$, we say $\lambda$-twistedly for
$\OO_\lambda$-twistedly.
 
\smallskip
\noindent
{\sl Remark.} An $\OO_G(\lambda =0)$-twistedly $G$-equivariant
$D$-module is just a $G$-equivariant $D$-module in [2, VII, 12.10], [24].
 
\vskip 10mm
\noindent
{\bf 2. Equivariant systems of linear partial differential equations}
 
\bigskip
All equivariant $D$-modules in these notes have the following forms:
 
\proclaim{Theorem}
Let $G$ acts on $X$ and $\lambda$ be a Lie algebra character. Let
$I$ be a finitely generated $G$-stable left ideal of the algebra of
linear differential operators $D_X(X)$. Then
$$
M = D_X/(I + \sum_{\theta \in \frak g} D_X (L_\theta - \lambda
(\theta)))
$$
is $\lambda$-twistedly $G$-equivariant ($L_\theta$ is the vector field on
$X$ given by the $G$-action:
$(L_\theta f)(x) = \displaystyle{\frac{d}{dt}}f(e^{-t \theta}x)|_{t=0}$).
\endproclaim
 
\demo{Proof}
Let $\OO_\lambda = D_G v_\lambda = \OO_G v_\lambda \  (L_{\theta ,G}
v_\lambda = \lambda (\theta ) v_\lambda , \  L_{\theta ,G}$ is the
right invariant vector field on $G$ corresponding to $\theta$)
and $u$ the generator of $M$, $M = D_X u \  (I \, u = (L_\theta -
\lambda (\theta ))u = 0$).
Take $v_\lambda \boxtimes u$ a generator of $\OO_\lambda \boxtimes M$
and $\tilde u = 1 \otimes u \in (\alpha^* M)(X)$.
 
First we see the actions of $L_{\theta ,G} \boxtimes 1$ and
$1 \boxtimes L_\theta$ on $\alpha^* M$. Taking a local coordinate
system $\{ x_i, \partial_i \}$ on $X$, we have
$\xi \tilde u = \sum_i \xi (x_i \circ \alpha ) \otimes \partial_i u$
for a vector field $\xi$ on $G \times X$. But then, by the right
invariance of $L_{\theta ,G}$,
$$
L_{\theta ,G} (x_i \circ \alpha ) = (L_\theta x_i) \circ \alpha.
$$
Hence
$$
\split
(L_{\theta ,G} \boxtimes 1) \tilde u
&= \sum_i 1 \otimes (L_\theta x_i) \partial_i u\\
&= 1 \otimes L_\theta u\\
&= \lambda (\theta )(1 \otimes u).
\endsplit
$$
 
For $1 \boxtimes L_\theta$, on the slice $g \times X \hookrightarrow X$,
$$
L_\theta (x_i \circ \alpha )|_{g \times X} = L_{g \theta} x_i
\quad (g \in G).
$$
Hence on $X \simeq g \times X \hookrightarrow G \times X$,
$$
\split
(1 \boxtimes L_\theta) \tilde u |_{g \times X}
&= \sum_i L_\theta (x_i \circ \alpha ) \otimes
\partial_i u|_{g \times X}\\
&= \sum_i L_{g \theta} (x_i) \partial_i u\\
&= L_{g \theta }u\\
&= \lambda (g \theta )u\\
&= \lambda (\theta )u
\endsplit
$$
where $g \theta = \roman{Ad}g \, \theta$ is the adjoint action of $g$
on $\frak g$ and the last equality comes from the $G$-invariance
of $\lambda \in {\frak g}^*$ (Lie algebra character).
This means that the section
$(1 \boxtimes L_\theta - \lambda (\theta )) \tilde u$
takes a value zero on every slice $g \times X$.
Since this section is in an $\OO_{G \times X}$-coherent subsheaf
$F_1D_{G \times X} \tilde u$,
$$
(1 \boxtimes L_\theta - \lambda (\theta )) \tilde u = 0
$$
by the Nakayama lemma.
 
Secondly for $P \in I$, we also have
$$
(1 \boxtimes P) \tilde u |_{g \times X} = P^g
\quad \text{on } g \times X
$$
by the same computation ($P^g$ is the $g$-translate of $P$ under
the $G$-action on $D_X(X)$). Since $I$ is $G$-stable, $P^g \in I$
and hence
$$
(1 \boxtimes P) \tilde u |_{g \times X} = 0
\quad \text{for } g \in G.
$$
Hence by the same reason as above,
$$
(1 \boxtimes P) \tilde u = 0 \quad \text{for } P \in I.
$$
 
We have thus proved $\roman{Ann} \, v_\lambda \boxtimes u \subset
\roman{Ann} \, \tilde u$ and hence the natural $D$-module
homomorphism
$$
\phi : \OO_\lambda \boxtimes M \longrightarrow \alpha^* M.
$$
By computation similar to the above, for the filtration
$F_mM = \mathbreak (F_mD_X)u \  (m \geq 0)$,
$$
\phi_m : \OO_\lambda \boxtimes F_mM \longrightarrow \alpha^*F_mM
$$
is surjective. Again restricting $\phi_m$ on the slice $g \times X$,
$$
\phi_m|_{g \times X} : \CC \boxtimes F_mM \  
\widetilde{\longrightarrow} \  (\alpha^*F_mM)_{g \times X}
\  \widetilde{\longrightarrow} \  F_mM.
$$
Thus again by the Nakayama lemma, $\phi_m$ is an isomorphism
and hence so is $\phi$. q.e.d.
\enddemo

 
 
 

\vskip 10mm
\noindent
{\bf 3. The Harish-Chandra system for characters}
 
\bigskip
Let $G$ be a connected reductive algebraic group over $\CC$ and
$Z$ be the center of the enveloping algebra $U(\frak g)$ of the
Lie algebra $\frak g = \roman{Lie} \, G$. Fix an algebra homomorphism
$$
\chi : Z \longrightarrow \CC.
$$
Consider the $D_G$-module $M_\chi = D_G u$ defined by
$$
M_\chi : \left\{
\alignedat 2 &(\partial_z - \chi (z))\, u = 0 &&\qquad (z \in Z) \\
&(L_\theta + R_\theta )\, u = 0 &&\qquad (\theta \in \frak g ),
\endalignedat \right.
$$
where $\partial_z$ is the two-sided invariant differential operator on
$G$ corresponding to $z \in Z$ and $L_\theta$ (resp. $R_\theta$) is
the right (resp. left) invariant vector field corresponding to
$\theta \in \frak g$. (Sorry! The symbols $L, R$ may be converse to the
usual ones. They follow the earlier notation $L_\theta =L_{\theta ,G},
\  R_\theta = R_{\theta ,G}$.) The vector field $L_\theta + R_\theta$
corresponds to the vector field arising from the inner action of $G$
on $G$ itself.
 
Since $Z$ is the subalgebra of $D_G(G)$ invariant under the two-sides
actions of $G$ and $L_\theta + R_\theta$ corresponds to the earlier
$L_\theta$ on $G$ by the inner action, $M_\chi$ is (0-twistedly)
$G$-equivariant under the inner action.
 
A distribution character of an irreducible admissible representation
of a real form $G_{\Bbb R}$ of $G$ is a solution to $M_\chi$ for some
$\chi$ (infinitesimal character). In order to analyze the behaviors
of characters, Harish-Chandra extensively investigated this system of
partial differential equations [10]
and we call it the Harish-Chandra system.
 
We shall see $M_\chi$ is a holonomic $D_G$-module. Since $G$ is an affine
variety, the $D_G$-module $M_\chi$ is the localization of the module of
global sections $M_\chi(G) = \Gamma (G, M_\chi )$. Put $D(G) = D_G(G)$
and
$$
I_\chi = \sum_{z \in Z}D(G)(\partial_z - \chi (z))
+ \sum_{\theta \in \frak g}D(G)(L_\theta + R_\theta ).
$$
Then
$$
M_\chi (G) = D(G)/ I_\chi.
$$
The filtration $F$ of $M_\chi (G)$ arising from the order filtration
of $D(G)$ is good and
$$
\roman{gr}^FM_\chi (G) = \roman{gr} \, D(G) / \roman{gr} \, I_\chi
$$
as is seen in I, 3.
 
Now let $\frak g$ be identified with the left invariant vector
fields ($\theta \leftrightarrow R_\theta$). Then under this
identification the cotangent bundle $T^*G$ is trivialized
as $G \times {\frak g}^*$ and
$$
\roman{gr} \, D(G) \simeq \CC [G] \otimes_{\CC} \CC [{\frak g}^*].
$$
On the other hand, by the Poincar\'e-Birkhoff-Witt theorem, we have
$$
\roman{gr} \, U(\frak g) \simeq S(\frak g) \simeq \CC [{\frak g}^*]
$$
where $S(\frak g)$ is the symmetric algebra over $\frak g$.
Considering the $G$-action on $U(\frak g)$ and $S(\frak g)$ arising
from the adjoint action, since the center $Z$ is the subalgebra of
$G$-invariants in $U(\frak g)$, we have
$$
\roman{gr} \, Z \simeq S(\frak g)^G \simeq \CC[{\frak g}^*]^G
$$
as the subalgebras of the above three algebras (complete reducibility
of $G$). Hence the symbols of $\partial_z - \chi (z)$ correspond
to the subset $\CC [{\frak g}^*]^G_+$ where $\CC [{\frak g}^*]^G_+ =
\CC [{\frak g}^*]^G \cap \CC [{\frak g}^*] \frak g$
in $\CC [{\frak g}^*]$.
 
On the other hand, the symbol of a vector field
$L_\theta + R_\theta$ at $(x, \xi ) \in G \times {\frak g}^*$
is written as
$$
\sigma_1(L_\theta + R_\theta )(x, \xi) = \langle - \roman{Ad}(x) \theta
+ \theta , \xi \rangle
$$
where $\langle \, , \, \rangle$ denotes the natural pairing
on $\frak g$ and $\frak g^*$.
Thus the ideal $\roman{gr} \, I_\chi$ contains $\CC [{\frak g}^*]^G_+$
and $\xi - x \xi$ on
$G \times {\frak g}^* \  (x \xi =
{}^t \roman{Ad}(x)\xi )$. Hence for the characteristic variety
$\roman{ch} \, M_\chi$ we have the inclusion relation
$$
\roman{ch} \, M_\chi \subset X = \{ (x, \xi ) \in G \times
{\frak g}^* \, | \, P(\xi)=0 \  (P \in \CC [{\frak g}^*]^G_+ ),
\, x \xi = \xi \, \}.
$$
 
But then by Kostant's theorem [18], the ideal generated by
$\CC [{\frak g}^*]^G_+$ in
$\CC [{\frak g}^*]$ is the defining
ideal of the nilpotent variety $N$ in ${\frak g}^*$ (the set of
nilpotent elements under the identification $\frak g \simeq
{\frak g}^*$ by a non-degenerate invariant bilinear form).
Thus we have
$$
X = \{ (x, \xi ) \in G \times N \, | \, \xi = x\xi \}.
$$
 
We now see an irreducible component of $X$ is of dimension equal to
$\dim G$. By Dynkin-Kostant, the nilpotent variety $N$ splits into
finitely many $G$-orbits under the coadjoint action:
$$
N = \coprod_{i=1}^{r}O_G(\xi_i).
$$
If $q : X \to N$ is the projection onto the second factor, then
$$
q^{-1}O_G(\xi_i) = \{ (x, \xi ) \in G \times O_G(\xi_i) \, |
\, x \in Z_G(\xi_i) \}
$$
is a fiber bundle over $O_G(\xi_i)$ with standard fiber
$Z_G(\xi_i) = \{x \in G | x \xi_i = \xi_i \}$.
Since $\dim O_G(\xi_i) = \dim G - \dim Z_G(\xi_i)$,
$\dim q^{-1}O_G(\xi_i) = \dim G$ and an irreducible component of $X$
is one of the closures of $q^{-1}O_G(\xi_i)$.
 
Thus we have seen $\dim \roman{ch} \, M_\chi \leq \dim G$ and hence
$M_\chi$ is holonomic.
 
More precisely, the following is known:
 
\proclaim{Theorem}
$$
\roman{ch}\, M_\chi = X.
$$
Furthermore the characteristic cycle of $M_\chi$ (for definition,
see [11, II]) is given by the intersection cycles
of $V$ and $G \times N$ where
$V = \{ (x, \xi ) \in G \times {\frak g}^* \, | \, x\xi = \xi \}$
is the commuting variety.
\endproclaim
 
\demo{Proof} See [14] for a Lie algebra version and [17]
for a group version.
\enddemo
 
Let $G_{\roman{rs}}$ be the set of regular semisimple elements $s \in G$,
i.e., semisimple and the centralizer $Z_G(s)$ of $s$ in $G$
has the dimension equal to
$\roman{rank} \, G$. Then for $s \in G_{\roman{rs}}$, if $\xi \in N$ and
$s\xi= \xi$, then $\xi=0$.
Thus $\roman{ch} \, M_\chi |_{G_{\roman{rs}}} \subset X \cap
(G_{\roman{rs}} \times {\frak g}^*) = G_{\roman{gr}} \times 0
= T_{G_{\roman{rs}}}^{*}(G_{\roman{rs}})$ and hence by I, 5.2,
$M_\chi|_{G_{\roman{rs}}}$ is a connection.
There is a well-known formula by
Harish-Chandra describing this connection on $G_{\roman{rs}}$.
We shall introduce this formula in the $D$-module language.
 
Fix a maximal torus $T$ in $G$ and let
$T_{\roman{reg}} = T \cap G_{\roman{rs}}$ the set of regular elements
in $T$.
Put $\widetilde{G_{\roman{rs}}} = G/T \times T_{\roman{reg}}$
and define the map
$$
p : \widetilde{G_{\roman{rs}}} \longrightarrow G_{\roman{rs}}
$$
by $p(gT, t) = gTg^{-1}$. Let $W=N_G(T)/T$ be the Weyl group for $T$.
Then $W$ acts on $\widetilde{G_{\roman{rs}}}$ by
$(gT, t)\dot w = (g \dot w T, w^{-1}tw) \  (w= \dot w T \in N_G(T))$
and $p$ is a Galois covering with group $W$
($\widetilde{G_{\roman{rs}}}/W \simeq G_{\roman{rs}}$).
 
We shall describe the inverse image $p^*(M_\chi |_{G_{\roman{rs}}})$ of
the connection $M_\chi |_{G_{\roman{rs}}}$.
Let $R_+$ be the set of positive roots and put the difference
$$
\Delta = \prod_{\alpha \in R_+}(t^{\alpha /2} - t^{- \alpha /2})
\quad (t \in T).
$$
($\Delta$ is possibly a two-valued function on $T$ and $T_{\roman{reg}}=
\{ t | \Delta (t) \not= 0 \}$.) Let $\frak t = \roman{Lie}\, T$ and
$S(\frak t)$ the symmetric algebra over $\frak t$ which is identified
both with invariant differential operators on $T$ and with the
polynomial algebra on ${\frak t}^*$. Denote by $S(\frak t)^W$ the
subalgebra of $W$-invariants in $S(\frak t)$.
 
For $\lambda \in {\frak t}^*$, consider the $D_{T_{\roman{reg}}}$-module
$M_\lambda^{\roman{rad}}$ defined by
$$
\Delta^{-1}(P(\partial ) - P(\lambda + \rho )) \Delta v =0
\quad (P \in S(\frak t)^W).
$$
where $\rho \in \frak t^*$ is the half sum of positive roots in $R_+$.
($P(\partial )$ is considered as an invariant differential operator
and $P(\lambda + \rho )$ is the value of the polynomial $P \in
\CC [\frak t^*]$ at $\lambda + \rho \in \frak t^*$.)
It is easily seen that $M_\lambda^{\roman{rad}}$ is a connection
on $T_{\roman{reg}}$.
 
By the standard theory of Chevalley and Harish-Chandra, the set
of algebra homomorphisms $\roman{Hom}_{\CC -\roman{alg}}(Z, \CC )$
is identified with the $W$-quotient ${\frak t}^*/.W$ where the
$.W$-action on $\frak t^*$ is defined by
$$
w.\lambda = w(\lambda + \rho ) - \rho
\quad (w \in W, \  \lambda \in \frak t^*).
$$
Hence for an infinitesimal character
$\chi \in \roman{Hom}_{\CC -\roman{alg}}(Z, \CC )$,
there exists $\lambda \in \frak t^*$ uniquely up to the $.W$-action.
 
\proclaim{Theorem (Harish-Chandra)}
Let $\chi_\lambda$ be the infinitesimal character corresponding to
$\lambda \in \frak t^*$ as above. Then we have a $D$-module
isomorphism on $G \times T_{\roman{reg}}$
$$
p^*(M_{\chi_\lambda}|G_{\roman{rs}}) \simeq
\OO_{G/T} \boxtimes M_\lambda^{\roman{rad}}.
$$
\endproclaim
 
It is rather easy to see that $M_\lambda^{\roman{rad}}$ and hence
$M_\chi |_{G_{\roman{rs}}}$ is a regular connection in the sense
of Deligne (I, 5.2). In the $D$-module theory, there exists a
unique minimal extension $(M_\chi |_{G_{\roman{rs}}})^\sim$ of
$M_\chi |_{G_{\roman{rs}}}$ on $G$ (related to the intersection
cohomology theory of Goresky-MacPherson) and this extension
$(M_\chi |_{G_{\roman{rs}}})^\sim$ is a regular $D_G$-module.
Finally, we close this section by citing the following:
 
\proclaim{Theorem ([14], [17])}
$$
M_\chi \simeq (M_\chi |_{G_{\roman{rs}}})^\sim.
$$
In particular, $M_\chi$ is a regular $D_G$-module.
\endproclaim
 
For applications of these considerations to the representation theory,
we refer to [12], [13], [14], [15].
 
\vskip 10mm
\noindent
{\bf 4. $D$-modules defined by a linear action and its orbits}
 
\bigskip
Let $G$ be an algebraic group acting linearly on a vector space $V$.
The dual space $V^*$ of $V$ is then acted by $G$ through the
contragredient action. Fix a special reference point $\xi \in V^*$ and
let $I_\xi = I(O_G(\xi)) \subset \CC [V^*]$ be the ideal consisting
of functions vanishing on $O_G(\xi)$, the $G$-orbit of $\xi$.
($I_\xi$ is the defining ideal of the orbit closure
$\overline{O_G(\xi)}$.)
Under the identification $\CC [V^*] \simeq S(V), \  I_\xi$ is
a $G$-stable ideal in the algebra $S(V)$ of linear diferential operators
with constant coefficients. Hence if we fix a Lie algebra character
$\lambda : \frak g \to \CC$, the following defines a $\lambda$-twistedly
$G$-equivariant $D_V$-module $M_{\lambda , \xi}$:
$$
M_{\lambda , \xi} : \left\{
\alignedat 2 &(L_\theta - \lambda (\theta ))\, u = 0
&& \qquad (\theta \in \frak g) \\
&P(\partial ) \, u = 0 && \qquad (P(\partial ) \in I_\xi ),
\endalignedat \right.
$$
 
\smallskip
\noindent
{\bf Example 1.}
Let $V$ be a prehomogeneous vector space under group $G$ (i.e., $V$
has a dense $G$-orbit). Let $\chi : G \to \CC^\times$ be an algebraic
group homomorphism and choose $\lambda \in \CC$.
For a Lie algebra homomorphism
$$
\lambda d \chi : \frak g \longrightarrow \CC ,
$$
$M_\lambda = D_Vu$ defined by
$$
(L_\theta - \lambda d \chi (\theta ))\, u =0 \quad (\theta \in \frak g)
$$
is $\lambda d\chi$-twistedly $G$-equivariant.
In many cases, $M_\lambda$ is known to be regular holonomic
and if $f$ is a relative invariant
for $\chi$ ($f(gx)= \chi (g)f(x)$), then $f^\lambda$ is a solution
to $M_\lambda$. The $D_V$-module $M_\lambda$ is important
for the study of the $b$-function of $f$. For recent results in this
field, refer to works of M. Muro and A. Gyoja ([21], [9]).
 
\smallskip
\noindent
{\bf Example 2.}
The Lie algebra version of the Harish-Chandra system for invariant
eigendistributions introduced in the previous section is simply an
example of this kind.
Let $G$ be a connected reductive algebraic group over $\CC$ and
$\frak g$ its Lie algebra as in the previous section 3.
Let $\xi$ be an regular element in $\frak g^*$ (which means that
the centralizer $Z_G(\xi)$ has the dimension equal to the rank of
$\frak g =\dim T$).
Consider the (0-twistedly) $G$-equivariant $D_{\frak g}$-module
$M_\xi = M_{0,\xi}$ corresponding to the regular orbit $O_G(\xi )$
on the Lie algebra $\frak g$.
 
Now by [18], it is known that
$I_\xi = I(\overline{O_G(\xi )})$ is generated by $P - P(\xi )
\  (P \in \CC [\frak g^*]^G)$. Hence $M_\xi$ is defined by
$$
M_\xi : \left\{
\alignedat 2 &(P(\partial ) - P(\xi ))\, u = 0 && \qquad (P(\partial ) \in
S(\frak g)^G) \\
&L_\theta u = 0 && \qquad (\theta \in \frak g),
\endalignedat \right.
$$
where we apply the identification
$$S(\frak g)^G \simeq \CC [\frak g^*]^G.$$
The above $D_{\frak g}$-module $M_\xi$ is of course defined for arbitrary
$\xi \in \frak g^*$ but even if $\xi$ is not regular, it gives the same
$D$-module by the following reason. If $\xi = \xi_s + \xi_n$ is the
Jordan decomposition of the element $\xi$, then $P(\xi ) = P(\xi_s)$
for any $P \in \CC [\frak g^*]^G$. Hence if $\xi_{\roman{reg}}$ is
a regular element such that the semisimple part of $\xi_{\roman{reg}}$
coincides with $\xi_s$ (such regular elements are unique up to the
$G$-action), then $M_\xi = M_{\xi_{\roman{reg}}}$.
 
This $D_{\frak g}$-module is investigated in details in [14].
As in 3, it is not difficult to see that this is holonomic but
in contrast with the group case in 3, $M_\xi$ is regular (in the
algebraic sense) only for a nilpotent $\xi$,
which also follows from the theorem in the next section.
 
\smallskip
\noindent
{\bf Example 3.}
Let $G =T$ be an algebraic torus ($\simeq {\CC}^{\times n}$) acting
linearly on a vector space $V$. Assume the action is faithful and $T$
contains all homotheties on $V$. Consider the contragredient action
of $T$ on $V^*$ and take a full diagonal subgroup $D \subset GL(V^*)$
such that $T \hookrightarrow D$. Take a reference point $\text{\bf \i}
 \in V^*$ such that
its $D$-orbit $D. \text{\bf \i} = O_D(\text{\bf \i})$ is open in $V^* 
\  (D \simeq
D. \text{\bf \i} \simeq {\CC}^{\times N} \  (\dim V = N))$. Then for
$\lambda \in \frak t^* \  (\frak t = \roman{Lie}\, T \simeq \CC^n )$,
$$
M_\lambda : \left\{
\alignedat 2 &(L_\theta - \lambda (\theta ))\, u = 0
&&\qquad (\theta \in \frak t) \\
&P \, u = 0 &&\qquad (P \in I(\overline{O_T(\text{\bf \i})}))
\endalignedat \right.
$$
gives a $\lambda$-twistedly $T$-equivariant $D_V$-module.
 
Gelfand [7] calls $M_\lambda$ a generalized hypergeometric
system of linear partial differential equations. In the final
section, we shall look at this $D$-module in more details.
In particular, this turns out to be a regular holonomic $D_V$-module.
 
Gelfand and his other colaborators also consider certain equivariant
holonomic $D$-modules which arise from unipotent group actions and
call them generalized Airy equations (irregular at $\infty$) [6].

 
 
 

\vskip 10mm
\noindent
{\bf 5. Some regularities}
 
\bigskip
 
In this section, we prove the regularity of a certain equivariant
$D$-module, which is suggested by a conversation with M.~ Kashiwara
whom we sincerely thank. Here we assume some standard
knowledge of algebraic $D$-modules as in [1], [2], [11].
 
Let $G$ be a connected algebraic group acting on a smooth algebraic
variety $X$.
Let $\lambda : \frak g \to \CC$ be a Lie algebra character and
$\OO_\lambda$ the rank one connection on $G$ defined by $\lambda$ in 2.
In case $\lambda = 0$, the following theorem is in [11, VII, 12.11].
 
\proclaim{Theorem}
Let $M$ be a $\lambda$-twistedly $G$-equivariant coherent $D_X$-module
whose support $\roman{Supp}\, M$ (support in the sheaf theory) consists
of finitely many $G$-orbits in $X$. If the twisting datum $\OO_\lambda$
is a regular connection on $G$, then $M$ is regular holonomic.
\endproclaim
 
\demo{Proof}
We use the induction on the number of $G$-orbits in $\roman{Supp}\, M$.
First we see the case that the orbit number is 1. Let
$$
G \overset \pi\to\longrightarrow \roman{Supp}\, M
= O_G(o)\overset i\to\hookrightarrow X
$$
($o \in \roman{Supp}\, M, \  i$ is a closed immersion).
Since by the Kashiwara lemma $M \simeq \int_i i^!M$
($i^! = \Cal H^0 i^!$ is the twisted inverse functor)
and $i^!M$ is again $G$-equivariant, it is enough to show
that $i^!M$ is regular holonomic. Thus we may assume that
$\roman{Supp}\, M=X=O_G(o)$.
Considering the maps
$$
G\overset \iota\to\simeq {G \times o}\overset{1 \times p}
\to\hookrightarrow {G \times X}\overset \alpha\to\longrightarrow X,
$$
we have
$$
\pi = \alpha \circ (1 \times p) \circ \iota
$$
where $p : o \hookrightarrow X$ is the inclusion map. Hence
$$
\align
\pi^*M &\simeq \iota^*(1 \times p)^* \alpha^*M \\
&\simeq \iota^*(1 \times p)^*(\OO_\lambda \boxtimes M) \\
&\simeq \OO_\lambda \boxtimes p^*M.
\endalign
$$
Since $M$ is $D$-coherent, $p^*M$ is a finite dimensional vector space
over $\CC$. Therefore $\pi^*M$ is a finite sum of a regular connection
$\OO_\lambda$. Since $\pi$ is a smooth map, $M$
is thus regular holonomic.
 
Secondly, we are going into a general case. Take a $G$-stable open set
$U$ in $X$ such that $U \cap \roman{Supp}\, M$ is the union of all
maximal dimensional $G$-orbits in $\roman{Supp}\, M$. Let
$$
U \overset j\to\hookrightarrow X \overset i\to\hookleftarrow Y
= X \setminus U
$$
be the corresponding open and closed immersions. Then we have the
exact triangle
$$
R\Gamma_Y(M) \longrightarrow M \longrightarrow j_*j^*M
\overset{+1}\to\longrightarrow R\Gamma_Y(M)
$$
in the derived category of coherent $D$-modules.
Note that all cohomology $D$-modules of the complexes $R\Gamma_Y(M)$
and $j_*j^*M$ are again $\lambda$-twistedly $G$-equivariant.
$\roman{Supp}\, j^*M$ is a disjoint union of $G$-orbits of the same
dimensions and hence by the case of orbit number = 1, $j^*M$ is
regular holonomic. Thus $j_*j^*M$ is a regular holonomic complex of
$D$-modules. On the other hand, the number of $G$-orbits in
$\roman{Supp}\, R\Gamma_Y(M)$ is less than that of $\roman{Supp}\, M$
and hence by the induction $R\Gamma_Y(M)$ is regular holonomic.
Thus the remaining vertex $M$ in the exact triangle is regular holonomic.
q.e.d.
\enddemo
 
Now we consider the case discussed in Section 4. Let $G$ act on a
vector space $V$ and take $\zeta \in V^*$ and a Lie algebra character
$\lambda$. Denote by $M_{\lambda ,\zeta}$ the
$D(V) = D_V(V)$(Weyl algebra)-module $D(V)\, u$ defined by
$$
\left\{
\alignedat 2
&(L_\theta - \lambda (\theta ))\, u = 0
&&\qquad (\theta \in \frak g) \\
&P(\partial )\, u = 0 &&\qquad (P(\partial ) \in I_\zeta )
\endalignedat \right.
$$
for $I_\zeta = I(\overline{O_G(\zeta )})$.
 
In general, for a Weyl algebra module $M$ on $V$, let $\hat M$ be
its Fourier transform. That is, $\hat M$ is a $D(V^*)$-module on
the dual space $V^*$ defined by the substitutions
$x_i \mapsto -\partial_{\xi_i}, \  \partial_{x_i} \mapsto
\xi_i$
for the dual coordinate systems $(x_i)$ in $V$ and
$(\xi_i)$ in $V^*$. The Fourier transform $\widehat{M_{\lambda ,\zeta}}$
is then defined by
$$
\left\{
\alignedat 2
&(\widehat{L_\theta} - \lambda (\theta ))\, \hat u = 0
&&\qquad (\theta \in \frak g) \\
&P(\xi )\, \hat u = 0 &&\qquad (P \in I_\zeta )
\endalignedat \right.
$$
where $\widehat{L_\theta}$ is the operator after the above
Fourier substitutions and $\hat u$ is the generator of
$\widehat{M_{\lambda ,\zeta}}$.
 
A Weyl algebra module $M$ is said to be {\it homogeneous} if the action
of the Euler operator $\sum_i x_i\partial_{x_i}$ is locally nilpotent
on $M$. If $M = D(V)/I$ and if the ideal $I$ is generated by homogeneous
elements in the gradation of $D(V)$ by $\roman{deg}\, \partial_{x_i}
= 1$ and $\roman{deg}\, x_i = -1$, then $M$ is homogeneous in this
sense.
 
It is known that a homogeneous Weyl algebra module $M$ is regular
holonomic if and only if its Fourier transform $\hat M$ is regular
holonomic (see for example [3, Th.7.4], [19]).
 
In the above examples, if $L_\theta$ has homogeneous degree 0 and the
ideal $I_\zeta$ is homogeneous
($\Leftrightarrow \  \overline{O_G(\zeta )}$
is conic), then $M_{\lambda ,\zeta}$ (and $\widehat{M_{\lambda ,\zeta}}$)
is homogeneous. Thus we have the following as a corollary of Theorem.
 
\proclaim{Corollary} Assume that $\OO_\lambda$ is a regular connection
(this is the case when $G$ is reductive). If $M_{\lambda ,\zeta}$ is
homogeneous and $\overline{O_G(\zeta )}$ consists of finitely many
$G$-orbits, then $M_{\lambda ,\zeta}$ is a regular holonomic $D_V$-module.
\endproclaim
 
\smallskip
\noindent
{\bf Examples.}
In Example 2 in Section 4, let $\xi$ be regular nilpotent. Then
$\overline{O_G(\xi )}$ is the nilpotent variety $N$ which is conic
and consists of finitely many orbits. Thus $M_\xi = M_0$ is regular
holonomic.
 
Example 3 (generalized hypergeometric systems) is also in this
category as will be seen in the next section.
 
\vskip 10mm
\noindent
{\bf 6. The Gelfand generalized hypergeometric systems}
 
\bigskip
\noindent
6.1 Description.
 
\smallskip
 
In this section, we shall look into more details of the equations
introduced in Example 3 in Section 4. Let $T = \CC^{\times n}$ be
an $n$-dimensional algebraic torus acting linearly on an $N$-dimentional
vector space $V$. Choose a coordinate system $(z_j)$ in $V$ such that
the $T$-action is diagonalized and denote by $\chi = (\chi_{ij})$
the integral ($n \times N$) matrix expressing the linear representation
$$
\chi : T=\CC^{\times n} \longrightarrow \CC^{\times N},
$$
i.e.,
$$
t.z_j = \prod_{i=1}^{n}t_i^{\chi_{ij}}z_j
\quad (t=(t_1, \cdots , t_n) \in T, \  1 \leq j \leq N).
$$
(We denoted by $\chi$ also the homomorphism given by $\chi$.)
We assume that $\chi$ is injective ($\Leftrightarrow \  
\roman{rank}\, \chi = n$) and that the image $\chi (T)$ contains all
homotheties, which means that there exists an integral vector
$c=(c_1, \cdots , c_n) \in {\Bbb Z}^n$ such that
$$
c \chi = (1, \cdots , 1) \quad \text{(all entries are 1)}
\eqno (\bigstar )
$$
(homogenuity condition).
 
In the Lie algebra $\frak t = \roman{Lie}\, T$, we choose the basis
$t_i \partial_{t_i} \  (1 \leq i \leq n)$. Then the action of
$t_i\partial_{t_i}$ on $V = \CC^N$ is given by
$$
\sum_{j=1}^N \chi_{ij} z_j\partial_{z_j}
\qquad \text{for }\  1 \leq i \leq n.
$$
Thus if we fix a Lie algebra character $\lambda \in \frak t^*
\simeq \CC^n \  (\lambda(t_i\partial_{t_i})= \lambda_i \in \CC \ 
(1 \leq i \leq n))$, then the equations
$(L_\theta - \lambda (\theta ))\, u=0 \  (\theta \in \frak t)$
correspond to
$$
(\sum_{j=1}^N \chi_{ij}z_j\partial_{z_j} - \lambda_i)\, u=0
\qquad (1 \leq i \leq n).
$$
 
Secondly in the dual space $V^* \simeq \CC^N$ (by the dual coordinate
system $(\xi_j)$ to $(z_j)$), fix a reference point
$$
\text{\bf \i}=(1, \cdots , 1) \in \CC^N \quad \text{(all entries are 1)}.
$$
The $T$-orbit of {\bf \i} is given by
$$
O_T(\text{\bf \i}) = \{ (t^{\chi_1}, \cdots , t^{\chi_N}) \  |
\  t \in T \  \}
$$
where $t^{\chi_j}= \prod_{i=1}^n t_i^{\chi_{ij}} \  (1 \leq j \leq N)$.
Define the sublattice $L$ in $\Bbb Z^N$ by
$$
L = \roman{Ker}\, \chi = \{ a = (a_1, \cdots , a_N) \in \Bbb Z^N \  |
\   \sum_{j=1}^N \chi_{ij} a_j =0 \  (1 \leq i \leq n) \}.
$$
Then we have
$$
O_T(\text{\bf \i})=\{ \xi = (\xi_j) \in \CC^{\times N} \  |\  
\xi^a = 1 \  (a \in L) \}
$$
where $\xi^a = \prod_{j=1}^N \xi_j^{a_j}$.
Thus the defining ideal $I(\overline{O_T(\text{\bf \i})})$ in
$\CC [\xi_1, \cdots , \xi_N ]$ is generated by
$$
\boxed{\xi}_a = \prod_{a_j >0} \xi_j^{a_j} - \prod_{a_j<0}\xi_j^{|a_j|}
\qquad (a \in L).
$$
 
Hence the generalized hypergeometric system $M_\lambda \  (\lambda
\in \CC^N)$ in Example 3 in 4 is realized by the following system
of linear partial differential equations:
$$
M_\lambda : \left\{
\alignedat 2
&(\theta_i - \lambda_i)\, u=0 &&\qquad (1 \leq i \leq n) \\
&\square_a u=0 &&\qquad (a \in L= \roman{Ker}\, \chi)
\endalignedat \right.
$$
where $\theta_i = \sum_{j=1}^N \chi_{ij} z_j \partial_{z_j}$ and
$\square_a = \prod_{a_j>0} \partial_{z_j}^{a_j} -
\prod_{a_j<0} \partial_{z_j}^{|a_j|}$ for $a \in L$.
 
Note that the homogenuity condition $(\bigstar )$ implies
$\sum_{j=1}^N a_j =0$ for $a \in L$ and hence the operator
$\square_a$ is of homogeneous degree $\sum_{a_j>0}a_j =
\sum_{a_j<0} |a_j|$. Also the vector fields $\theta_i$ is of
homogeneous degree 0 by $\roman{deg}\, z_j =-1, \ 
\roman{deg}\, \partial_{z_j}=1$. Hence $M_\lambda$ is a homogeneous
Weyl algebra module in the sense of Section 5.
 
\bigskip
\noindent
6.2 Characteristic variety.
 
\smallskip
 
Let $D = D(V) =\CC [z_1, \cdots , z_N, \partial_{z_1}, \cdots ,
\partial_{z_N}]$ be the Weyl algebra over $\CC^N$ and write
$M_\lambda = D\, u$, the global section of the $M_\lambda$ in 6.1.
Then $I_\lambda = \roman{Ann}_D u$ is generated by $\theta_i-\lambda_i \ 
(1 \leq i \leq n)$ and $\square_a \  (a \in L)$. Taking the
gradation by the good filtration $F$ on $M_\lambda$ given by
the order filtration of $D$, we have
$$
\roman{gr}^F M_\lambda = \roman{gr}\, D /\roman{gr}\, I_\lambda
$$
where $\roman{gr}\, D = \CC [z_1, \cdots , z_N, \xi_1, \cdots , \xi_N]
= \CC [z, \xi]$.
Since $\roman{gr}\, I_\lambda$ contains
$$
\alignedat 2
&\bar \theta_i = \sum_{j=1}^N \chi_{ij} z_j \xi_j
&&\qquad (1 \leq i \leq n), \\
&\boxed{\xi}_a = \prod_{a_j>0} \xi_j^{a_j} - \prod_{a_j<0} \xi_j^{|a_j|}
&&\qquad (a \in L),
\endalignedat
$$
there exists a surjective homomorphism
$$
A \twoheadrightarrow \roman{gr}^F M_\lambda
$$
where $A=\CC [z, \xi]/(\bar \theta_i , \boxed{\xi}_a | 1 \leq i \leq n,
a \in L )$. Hence for the characteristic variety of $M_\lambda$
$$
\roman{ch}\, M_\lambda \subset \roman{Supp}\, A = \{ (z, \xi) \in
\CC^{2N} \  | \  \bar \theta_i =0, \boxed{\xi}_a =0 \ 
(1 \leq i \leq n, \  a \in L) \}.
$$
 
So far, we have regarded $\CC^{2N}$ as the cotangent bundle $T^*(V)$
of $V$ ($(\xi_j)$ is the fiber coordinate), but this same space can
also be regarded as the cotangent bundle $T^*V^*$ of the dual vector
space $V^*$ ($(z_j)$ is in turn the fiber coordinate). Now the orbit
closure
$$
\overline{O_T(\text{\bf \i})} = \{ \xi \in V^* \  |\  \boxed{\xi}_a
=0 \  (a \in L) \}
$$
in $V^*$ splits into finitely many $T$-orbits:
$$
\overline{O_T(\text{\bf \i})} = \coprod_{k=1}^r O_T(p_k).
$$
Then the equations $\bar \theta_i =0$ gives the conormal condition
at $\xi \in O_T(p_k)$ and hence we have the following.
 
\proclaim{Lemma}
$$
\roman{Supp}\, A = \coprod_{k=1}^r T^*_{O_T(p_k)}(V^*)
$$
where $T^*_X(V^*)$ denotes the conormal bundle of $X \subset V^*$.
In particular, every irreducible component of $\roman{ch}\, M_\lambda$
is the closure of the conormal bundle of a $T$-orbit $O_T(p_k)$
and hence $N$-dimensional.
\endproclaim
 
We thus have:
\proclaim{Theorem}
$M_\lambda$ is a regular holonomic $D_V$-module.
\endproclaim
 
\demo{Proof} The regularity follows from Corollary in 4, since
$M_\lambda$ is homogeneous and $\overline{O_T(\text{\bf \i})}$
consists of finitely many $T$-orbits.
\enddemo
 
\noindent
{\sl Remark.}
In [7], the authors claim the isomorphism $A \simeq \roman{gr}^F
M_\lambda$. Then the characteristic cycle of $M_\lambda$ coincides
with that of the commutative algebra $A$, which has a nice
description coming from the combinatorics of polytopes defined by
the integral column vectors in the matrix $\chi$.
 
We understand this isomorphism only in case that the
orbit closure $\overline{O_T(\text{\bf \i})}$ is a normal variety.
Most important examples seem to satisfy this normality condition.
In particular, Mutsumi Saito has checked this normality condition
for all systems arising from the symmetric pairs [22].
 
\bigskip
\noindent
6.3 Classical cases.
 
\smallskip
 
Kinds of classical hypergeometric functions are defined on the
quotient space (or its compactification) by the torus action,
instead on $\CC^N$. We shall look at this situation on the
generic part $\CC^{\times N}$. The torus action defines the injective
homomorphism of tori $\chi : T =\CC^{\times n} \hookrightarrow
\CC^{\times N}$. Let
$$
e \longrightarrow \CC^{\times n} \overset \chi\to\longrightarrow
\CC^{\times N} \overset \pi\to\longrightarrow \CC^{\times l}
\longrightarrow e
$$
be the exact sequence of tori ($\CC^{\times l}$ is the quotient and
$l= N-n$). The quotient map $\pi$ is also given by an integral
($N \times l$) matrix (denoted also by the same symbol $\pi$) if
a coordinate system on $\CC^{\times l}$ is fixed. We want a system
of differential equations on $\CC^l$ whose inverse image on
$\CC^{\times N}$ is related to Gelfand's $M_\lambda|_{\CC^{\times n}}$.
 
We take a heuristic view point. For a function $v$ on $\CC^{\times l}$,
set the function $u$ on $\CC^{\times N}$ as
$$
u(z) = z^\Lambda v(\pi(z)) \qquad
\text{for some }\  \Lambda = (\Lambda_1, \cdots , \Lambda_N) \in \CC^N
$$
where $z^\Lambda = \prod_{j=1}^N z_j^{\Lambda_j}$.
Choose a coordinate system $x=(x_k)$ in $\CC^l \  (1 \leq k \leq l)$
and let
$$
x_k =\pi (z)_k = \prod_{j=1}^N z_j^{\pi_{jk}} \qquad (1 \leq k \leq l)
$$
(here $\pi = (\pi_{jk})$ is an integral $N \times l$ matrix).
Put $\vartheta_{z_j} = z_j \partial_{z_j}$ and $\vartheta_{x_k} =
x_k \partial_{x_k}$. Then
$$
\partial_{z_j}u = z_j^{-1} z^\Lambda (\sum_{k=1}^l \pi_{jk}
\vartheta_{x_k} + \Lambda_j )\, v \qquad (1 \leq j \leq N)
$$
and hence
$$
\vartheta_{z_j} u = z^\Lambda (\sum_{k=1}^l \pi_{jk} \vartheta_{x_k}
+ \Lambda_j )\, v.
$$
Thus, if $u$ satisfies
$$
(\theta_i - \lambda_i )\, u = 0 \qquad (1 \leq i \leq n),
$$
then $v$ satisfies
$$
z^\Lambda (\sum_{j,k} \chi_{ij} \pi_{jk} \vartheta_{x_k} +
\sum_{j=0}^N \chi_{ij} \Lambda_j - \lambda_i )\, v =0.
$$
But then since $\chi \pi =0$ as a multiplication of matrices (by the
definition of $\pi$), in order to get a non-trivial solution for $v$
to the above equations, we have the linear matrix equation for
$\Lambda \in \CC^N$:
$$
\chi \Lambda = \lambda
$$
where the matrix $\chi$ and the vector $\lambda \in \CC^n$ are given.
 
In order to examine the second equations $\square_a u =0$, we
compute the following iterated differentiations:
$$
\partial_{z_j}^m u = z_j^{-m} z^\Lambda (D_j + \Lambda_j - m+1)_m v
$$
where $D_j = \sum_{k=1}^l \pi_{jk} \vartheta_{x_k}$ and
$$
(T)_m = T(T+1) \cdots (T+m-1)
\qquad (m \geq 1, \  T \text{: indeterminate}).
$$
Hence $\square_a u =0$ corresponds to
$$
\{ \prod_{a_j>0}(D_j +\Lambda_j -a_j +1)_{a_j} -
z^a \prod_{a_j<0}(D_j + \Lambda_j +a_j +1)_{|a_j|} \}\, v =0
\  (a \in L) \eqno (\sharp ).
$$
Since $a \in L = \roman{Ker}\, \chi = \roman{Im}\, \pi$,
if $\pi_k \in \Bbb Z^N$ is the $k$-th column vector of the $N \times l$
matrix $\pi \  (1 \leq k \leq l)$, then
$$
a = \sum_{k=1}^l m_k \pi_k \qquad \text{for some } m_k \in \Bbb Z.
$$
Hence
$$
z^a = \prod_{j=1}^N z_j^{a_j} = \prod_{k=1}^l x_k^{m_k}
\qquad (x_k = z^{\pi_k})
$$
and the equations $(\sharp )$ are defined on $\CC^l$. We thus have
the following:
 
\proclaim{Proposition}
For $\lambda \in \CC^n$, choose $\Lambda \in \CC^N$ such that
$\chi \Lambda = \lambda$. Let $N_\Lambda =D_{\CC^l} v$ on $\CC^l$
defined by the equations $(\sharp )$ for $a \in L$. Then the $D$-module
$\OO_\Lambda \otimes \pi^* N_\Lambda$ on $\CC^{\times N}$ is a quotient
$D$-module of $M_\lambda |_{\CC^{\times N}}$. ($\OO_\Lambda = D_{\CC^N}
z^\Lambda$).
\endproclaim
 
\smallskip
\noindent
{\sl Remark.}
Assume $l=N-n=1$ and let $\pi =(\pi_1 , \cdots ,\pi_N) \in \Bbb Z^N$.
Then the equations $(\sharp )$ are reduced to a single equation
for this $\pi$:
$$
\{\prod_{\pi_j>0}(\pi_j \vartheta_x + \Lambda_j - \pi_j +1)_{\pi_j}
- x \prod_{\pi_j<0}(\pi_j \vartheta_x +
\Lambda_j + \pi_j +1)_{|\pi_j|} \}\, v =0
$$
This ordinary differential equation is of Fuchsian type (regular in
our language) and has a solution of the ``classical generalized"
hypergeometric function of one variable ${}_pF_{p-1} \ 
(p= \sum_{\pi_j>0} \pi_j)$ for suitable parameters expressed by
$\pi$ and $\Lambda$.
 
For the further references in several variables ($l > 1$), see [7], [8].
 
 


\enddocument